# 不含3-圈的1-平面图的边染色

张欣,刘桂真,吴建良

(山东大学数学学院,山东 济南 250100)

**摘要**:利用权转移方法证明了最大度 $\Delta \geq 7$ 且不含3-圈的1-平面图是 $\Delta$-边可染的。
**关键词**:1-平面图;三角形;边色数
**中图分类号**:O157.5    **文献标志码**:A

## Edge coloring of triangle-free 1-planar graphs

ZHANG Xin, LIU Gui-zhen, WU Jian-liang

(School of Mathematics, Shandong University, Jinan 250100, Shandong, China)

**Abstract**: It is shown that each triangle-free 1-planar graph with maximum degree $\Delta \geq 7$ can be $\Delta$-colorable by the Discharging Method.

## 1 主要定理及其证明

本文仅考虑简单的有限无向图。设 $G$ 是一个图,分别用 $V(G), E(G), \Delta(G), \delta(G)$ 来表示它的点集合,边集合,最大度和最小度。特别地,若 $G$ 是一个平面图,则用 $F(G)$ 来表示它的面集合。其余定义若无特殊说明,均参考文献[1]。

一个图 $G$ 称为可 $k$-边可染的,当且仅当它的边集合可以用 $k$ 种颜色染色使得相邻的边染不同的颜色。图 $G$ 的边色数 $\chi'(G)$ 是使得 $G$ 是 $k$-边可染的最小 $k$ 值。对于边染色,Vising 和 Gupta 独立地证明了:对于简单图 $G$,均满足 $\Delta(G) \leq \chi'(G) \leq \Delta(G) + 1$。对于图 $G$,若 $\chi'(G) = \Delta(G)$,则称 $G$ 是第一类的;若 $\chi'(G) = \Delta(G) + 1$,则称 $G$ 是第二类的。若图 $G$ 是第二类的且对于 $G$ 中的任何一条边 $e$ 均满足 $\chi'(G - e) < \chi'(G)$,则称 $G$ 是临界的。如若 $G$ 是临界的且 $G$ 的最大度为 $\Delta$,则称 $G$ 是 $\Delta$-临界的。

一个图称为是1-平面的,当且仅当它可以画在一个平面上,使得它的任何一条边最多交叉另外一条边。与平面图相比,1-平面图的研究成果非常少。Ringle 猜想任何一个1-平面图均是可6-点可染的,此猜想已由 Borodin 在文献[2]中证实且猜想中的6是最优的。在文献[3]中,Borodin 还证明了任何一个1-平面图均是20-无圈点可染的。对于其它的染色,甚至是点染色的其它问题,均无任何研究。关于1-平面图的结构,Fabrici 在文献[4]中证明了:任意1-平面图均满足 $e(G) \leq 4v(G) - 8$;对于不含3-圈的1-平面图 $G$,有 $e(G) \leq 3v(G) - 6$。后者这个上界恰好与平面图的边数上界吻合。而在文献[5]中,Sanders 和 Zhao 证明了最大度 $\Delta \geq 7$ 的平面图是第一类的。本文将证明最大度 $\Delta \geq 7$ 且不含3-圈的1-平面图也是第一类的。

对于1-平面图 $G$,本文总是假设 $G$ 满足最优性,即已经将 $G$ 画在平面上使得 $G$ 中的交叉点的数目是最小的。因此对于交叉于 $z$ 点的两条边 $x_1y_1, x_2y_2$,它们的四个端点是两两不同的。设 $C(G)$ 是 $G$ 中所有交叉

---

点的集合(注意交叉点并非是 $G$ 中的顶点), $E_0(G)$ 是 $G$ 中没有被交叉的边的集合。定义 1-平面图 $G$ 的关联平面图 $G^\times$,使得 $V(G^\times) = V(G) \cup C(G), E(G^\times) = E_0(G) \cup \{xz, yz | xy \in E(G) \setminus E_0(G)\}$,其中 $z$ 是位于 $xy$ 上的交叉点。于是 $G$ 中的交叉点成为了 $G^\times$ 中的度数为 4 的顶点。在不引起混淆的情况下,我们仍然记 $G^\times$ 中的新点为交叉点,$G$ 中含有交叉点的边称为交叉边。由 1-平面图的定义知,$G^\times$ 中任何两个交叉点不相邻。若 $v$ 不是交叉点,则有 $d_{G^\times}(v) = d_G(v)$。因此下文若无特殊说明,对于非交叉点将不再区分 $d_{G^\times}(v)$ 与 $d_G(v)$,而是简单地用 $d(v)$ 来表示。下面将给出 $G^\times$ 的一个结构性质。对于 $v \in V(G^\times)$,用 $f_k(v)$ 表示 $v$ 在 $G^\times$ 中所关联的 $k$-面的数目,用 $n_c(v)$ 表示 $v$ 在 $G^\times$ 中所相邻的交叉点的个数。

**引理 1**　设 $G$ 是不含 3-圈的 1-平面图,则对于任意的 $v \in V(G)$ 均满足 $f_3(v) \leq n_c(v)$。

**证明**　由于 $G$ 是不含 3-圈的,从而 $v$ 所关联的 3-面必包含一个交叉点。若 $n_c(v) = 0$,则必有 $f_3(v) = 0$。若 $n_c(v) \neq 0$,取 $u$ 为 $G^\times$ 中与 $v$ 相邻的交叉点,则 $uv$ 最多关联 1 个 3-面。否则 $uv$ 关联 2 个 3-面,记为 $uvx$, $uvy$,从而 $vx, vy \in E(G)$。再由 $G$ 的最优性知,必有 $xy \in E(G)$,且 $u$ 是位于 $xy$ 上的交叉点。此时 $vxy$ 构成了 $G$ 中的一个 3-圈,矛盾。因此,与 $v$ 相邻的每一个交叉点最多关联一个包含 $v$ 的 3-面。由此即得 $f_3(v) \leq n_c(v)$。证毕。

**引理 2**　设 $G$ 是不含 3-圈的 1-平面图,对于 $v \in V(G)$,有如下结论:(1) 若 $d(v) = 2$,则 $f_3(v) = 0$;(2) 若 $d(v) = 3$ 且 $n_c(v) \geq 2$,则 $f_3(v) \leq 2$,其中等号成立当且仅当 $v$ 关联一个 $6^+$-面;(3) 若 $d(v) = 4$ 且 $n_c(v) \geq 3$,则 $f_3(v) \leq 2$;(4) 若 $d(v) = 5$ 且 $n_c(v) \geq 4$,则 $f_3(v) \leq 2$;(5) 若 $d(v) = 6$ 且 $n_c(v) \geq 5$,则 $f_3(v) \leq 2$;(6) 若 $d(v) = 7$ 且 $n_c(v) \geq 5$,则 $f_3(v) \leq 4$。

**证明**　这里仅详细给出(1)与(2)的证明,其余的可用与(2)同样的方法加以证明。假设 $d(v) = 2$ 且 $v$ 关联一个 3-面 $uvw$。则 $u, w$ 中必然有且仅有一个是交叉点,不妨为 $u$。故可设 $vv_1, ww_1$ 均为 $G$ 中的交叉边且二者在 $G$ 中交叉于点 $u$。现重新将边 $ww_1$ 画在与 $v, w, w_1$ 的面内使得 $ww_1$ 不再与路 $wvv_1$ 交叉。按如此方法可使得 $G$ 中的交叉点减少一个,此与 $G$ 的最优性矛盾。下证(2),当 $d(v) = 3$ 且 $n_c(v) \geq 2$ 时,记 $N_{G^\times}(v) = \{x, y, z\}$。先假设 $n_c(v) = 3$。此时必有 $f_3(v) = 0$,否则将出现两个交叉点相邻的情形,矛盾。当 $n_c(v) = 2$ 时,由引理 1 知 $f_3(v) \leq n_c(v) = 2$。若 $f_3(v) = 2$,不妨设 $xy, xz \in E(G^\times)$,此时只可能 $y, z$ 是与 $v$ 相邻的两个交叉点。设 $xx_1, xx_2$ 为 $G$ 中的交叉边,其中 $y, z$ 分别为这两条边上的交叉点。于是有 $x_1y, x_2z \in E(G^\times)$。由于 $G$ 是一个简单图,则 $x_1 \neq x_2$,否则将出现重边。若 $x_1x_2 \in E(G^\times)$,则 $xx_1xx_2$ 构成了 $G$ 中的一个 3-圈,从而 $x_1x_2 \notin E(G^\times)$。故 $v$ 关联的第三个面必为一个 $6^+$-面 $f$。证毕。

**引理 3**[6]　设 $G$ 是一个 $\Delta$-临界图, $uv \in E(G)$ 且 $d(v) = k \leq \Delta$,则 $w$ 相邻至少 $\Delta - k + 1$ 个 $\Delta$-度点。

**引理 4**[7]　设 $G$ 是一个 $k$-临界图,其中 $k \geq 8$,则 $e(G) \geq 3v(G)$。

**定理 1**　设 $G$ 是不含 3-圈的 1-平面图,若 $\Delta(G) \geq 7$,则 $\chi'(G) = \Delta(G)$,即 $G$ 是第一类的。

**证明**　假设定理不成立,则存在一个 $\Delta(G)$-临界的不含 3-圈的 1-平面图 $G$,从而有 $e(G) \leq 3v(G) - 6$。由引理 4 知,$\Delta(G) \leq 7$。而定理假设 $\Delta(G) \geq 7$,故有 $\Delta(G) = 7$。考虑图 $G^\times$,由欧拉公式知 $\sum_{v \in V(G^\times)} (d(f) - 4) + \sum_{f \in F(G^\times)} (d(f) - 4) = -8$。对于 $x \in V(G^\times) \cup F(G^\times)$,定义权如下:对于任意的 $v \in V(G^\times)$,令 $ch(v) = d(v) - 4$;对于任意的 $f \in F(G^\times)$,令 $ch(f) = d(f) - 4$。则 $\sum_{x \in V(G^\times) \cup F(G^\times)} ch(x) < 0$。下面对 $x \in V(G^\times) \cup F(G^\times)$ 中各元素的权值进行调整,记调整后的权值为 $ch'(x)$。由于权只是在各元素内部进行转移,从而权的加和保持不变,即 $\sum_{x \in V(G^\times) \cup F(G^\times)} ch'(x) = \sum_{x \in V(G^\times) \cup F(G^\times)} ch(x) < 0$。下面定义权转移规则:

R1　若 $uv \in E(G), d(u) = 7, 2 \leq d(v) \leq 6$,则 $u$ 向 $v$ 转移权值 $\dfrac{1}{d(v) - 1}$;

R2　若 $uv \in E(G), d(u) = 6$,当 $d(v) = 3$ 时,$u$ 向 $v$ 转移权值 $\dfrac{1}{2}$,当 $d(v) = 4$ 时,则 $u$ 向 $v$ 转移权值 $\dfrac{1}{6}$;

R3　若 $f = uvw$ 是 $G^\times$ 中的 3-面且 $w$ 为交叉点,则 $u, v$ 分别向 $f$ 转移权值 $\dfrac{1}{2}$;

R4　若 $f$ 是 $G^\times$ 中的一个 $6^+$-面,则 $f$ 向其关联的每个 3 度点转移权值 $\dfrac{2}{3}$。

下面计算按照如上转移规则后 $x \in V(G^\times) \cup F(G^\times)$ 中各元素的权值。设 $f \in F(G^\times)$,若 $d(f) = 3$,则 $f$

必关联一个交叉点,两个非交叉点,从而由 R3 知 $ch'(f) = ch(f) + 2 \times \frac{1}{2} = 0$。若 $4 \leq d(f) \leq 5$,则由以上五条规则知 $f$ 不转出也不转入权值,从而 $ch'(f) = ch(f) \geq 0$。若 $d(f) \geq 6$,由引理 3 知任何两个 3 度点均不相邻,从而 $f$ 最多关联 $\lfloor \frac{d(f)}{2} \rfloor$ 个 3 度点,故由 R4 得到, $ch'(f) \geq ch(f) - \frac{2}{3} \times \lfloor \frac{d(f)}{2} \rfloor \geq \frac{2d(f)-12}{3} \geq 0$。

设 $v \in V(G^\times)$,由于 $G$ 是临界的,显然有 $d(v) \geq 2$。若 $d(v) = 2$,从引理 3 知,$v$ 的两个邻点均为 7 度点,从而由 R1 知 $v$ 从每个邻点得到权值 1。由引理 2 中的(1)以及以上四条规则,$v$ 不向外转移权值。故有 $ch'(v) = ch(v) + 2 = 0$。若 $d(v) = 3$,由引理 3 知,$v$ 的邻点均为 6 度点或 7 度点。故由 R1 和 R2 知,$v$ 可从其邻点得到权值 $3 \times \frac{1}{2} = \frac{3}{2}$。由引理 2 中的(2)知 $v$ 最多关联两个 3-面。如果 $f_3(v) \leq 1$,则由 R4 知 $v$ 向外转移的权最多为 $\frac{1}{2}$,因此 $ch'(v) \geq ch(v) + \frac{3}{2} - \frac{1}{2} = 0$。如果 $f_3(v) = 2$,则由 R4 知 $v$ 向外转移的权最多为 $2 \times \frac{1}{2} = 1$。但此时 $v$ 还关联一个 $6^+$ 面 $f$,由 R4 知 $f$ 将向 $v$ 转移权值 $\frac{2}{3}$。因此 $ch'(v) \geq ch(v) + \frac{3}{2} - 1 + \frac{2}{3} > 0$。若 $d(v) = 4$,由引理 3 知,$v$ 的邻点均为 $5^+$-点。并且,若 $v$ 有一个 5 度邻点,则必然还有 3 个 7 度邻点;若有一个 6 度邻点,则另外 3 个邻点要么是 3 个 7 度点,要么是 2 个 7 度点和一个 6 度点。不论哪种情况,根据 R1 和 R2,$v$ 都将从其邻点得到权值 1。另外由引理 1 和引理 2 的(3)知 $f_3(v) \leq 2$,从而根据 R3,$v$ 最多向外转移权值 1,于是有 $ch'(v) \geq ch(v) + 1 - 1 = 0$。若 $d(v) = 5$,由引理 3 知 $v$ 必相邻至少 2 个 7 度点,从而根据 R1,$v$ 从其邻点得到的权值至少为 $\frac{1}{2}$。由引理 1 和引理 2 的(4)知 $f_3(v) \leq 3$,从而根据 R3,$v$ 最多向外转移权值 $\frac{3}{2}$,于是有 $ch'(v) \geq ch(v) + \frac{1}{2} - \frac{3}{2} = 0$。若 $d(v) = 6$,当 $v$ 与 1 个 3 度点或 4 度点相邻时,则由引理 3 知 $v$ 还相邻 4 个 7 度点,从而根据 R1 和 R2,$v$ 从其邻点中总共转入权值至少为 $4 \times \frac{1}{5} = \frac{4}{5}$,转出权值至多为 $\max\{\frac{1}{2}, 2 \times \frac{1}{6}\} = \frac{1}{2}$。当 $v$ 不相邻 3 度点和 4 度点时,由以上四条规则知 $v$ 既不从其邻点转入权值也不转出权值。另外,由引理 1 和引理 2 的(5)知 $f_3(v) \leq 4$,从而由 R3 知 $ch'(v) \geq ch(v) - 4 \times \frac{1}{2} + \min\{\frac{4}{5} - \frac{1}{2}, 0\} = 0$。若 $d(v) = 7$,则由引理 3 以及 R1 知 $v$ 向其邻点转移权值最多为 $\max\{1, 2 \times \frac{1}{2}, 3 \times \frac{1}{3}, 4 \times \frac{1}{4}, 5 \times \frac{1}{5}\} = 1$。而由引理 1 和引理 2 的(6)知 $f_3(v) \leq 4$,从而由 R3 知 $ch'(v) \geq ch(v) - 1 - 4 \times \frac{1}{2} = 0$。至此已经证得 $\sum_{x \in V(G^\times) \cup F(G^\times)} ch'(x) \geq 0$,矛盾。因此定理 5 成立,证毕。

**注**　文献[4]指出,存在一个含 3-圈的 7-正则 1-平面图 $G$,从而有 $\chi'(G) = \Delta(G) + 1$。因此,定理 1 中的不含 3-圈这个条件是必要的。